\documentclass[11pt]{article}
\usepackage{epic,latexsym,amssymb}
\usepackage{color}
\usepackage{tikz}

\newtheorem{thm}{Theorem}

\newtheorem{defn}{Definition}
\newtheorem{cor}{Corollary}
\newtheorem{obser}{Observation}
\newtheorem{prop}{Proposition}

\newcommand{\qed}{$\Box$}

\newcommand{\cL}{{\cal L}}
\newcommand{\cR}{{\cal R}}

\newcommand{\odd}{{\rm odd}}
\newcommand{\even}{{\rm even}}

\newcommand{\barD}{\overline{D}}
\newcommand{\bard}{\overline{d}}
\newcommand{\bardelta}{\overline{\delta}}
\newcommand{\barDelta}{\overline{\Delta}}
\newcommand{\barC}{\overline{C}}
\newcommand{\barG}{\overline{G}}
\newcommand{\barS}{\overline{S}}
\newcommand{\barv}{\overline{v}}

\newcommand{\modo}{{\rm mod}}

\newcommand{\oneD}{\gamma_{\times 1}}
\newcommand{\oneT}{\gamma_{\times 1,t}}
\newcommand{\kDom}{\gamma_{\times k}}
\newcommand{\kRDom}{\gamma_{\times k}^{r}}
\newcommand{\kTDom}{\gamma_{\times k,t}}
\newcommand{\kTRDom}{\gamma_{\times k,t}^{r}}
\newcommand{\oneR}{\gamma_{\times 1}^{r}}

\newcommand{\gt}{\gamma_t}
\newcommand{\gr}{\gamma_r}

\newcommand{\ktdom}{d_{\times k,t}}
\newcommand{\kdom}{d_{\times k}}
\newcommand{\kdoms}{d_{\times k}^*}
\newcommand{\krdom}{d_{\times k}^{\, r}}
\newcommand{\krtdom}{d_{\times k,t}^{\, r}}

\newcommand{\proof}{\noindent\textbf{Proof. }}
\newcommand{\2}{ \vspace{0.2cm} }
\newcommand{\1}{ \vspace{0.1cm} }

\let\oldenumerate\enumerate
\renewcommand{\enumerate}{
  \oldenumerate
  \setlength{\itemsep}{0pt}
  \setlength{\parskip}{0pt}
  \setlength{\parsep}{0pt}
}

\begin{document}

\title{$k$-Tuple Restrained Domination in Graphs}

\author{$^1$Michael A. Henning\thanks{Research
supported in part by the University of Johannesburg and the South African National Research Foundation.} \, and \, $^{2}$Adel P. Kazemi  \\ \\
$^{1}$Department of Pure and Applied Mathematics \\
University of Johannesburg \\
Auckland Park, 2006 South Africa\\
\small \tt Email: mahenning@uj.ac.za \\
\\
$^{2}$Department of Mathematics\\ University of Mohaghegh Ardabili \\ P.O.\ Box 5619911367, Ardabil, Iran\\
\small \tt Email: adelpkazemi@yahoo.com}

\date{}
\maketitle

\begin{abstract}
For $k \ge 1$ an integer, a set $S$ of vertices in a graph $G$ with minimum degree at least~$k-1$ is a $k$-tuple dominating set of $G$ if every vertex of $S$ is adjacent to at least $k-1$ vertices in $S$ and every vertex of $V(G) \setminus S$ is adjacent to at least $k$ vertices in $S$; that is, $|N_G[v] \cap S| \ge k$ for every vertex $v$ of $G$ where $N_G[v]$ denotes the closed neighborhood of $v$ which consists of $v$ and all neighbors of $v$. A $k$-tuple restrained dominating set of $G$ is a $k$-tuple dominating set $S$ of $G$ with the additional property that every vertex outside $S$ has at least $k$ neighbors outside $S$. The minimum cardinality of a $k$-tuple restrained dominating set of $G$ is the $k$-tuple restrained domination number of $G$. When $k=1$, the $k$-tuple restrained domination number is the well-studied restrained domination number. In this paper, we determine the $k$-tuple restrained domination number of several classes of graphs. Tight upper bounds on the $k$-tuple restrained domination number of a general graph are established. We present basic properties of the $k$-tuple restrained domatic number of a graph which is the maximum number of the classes of a partition of $V(G)$ into $k$-tuple restrained dominating sets of $G$.
\end{abstract}


\newpage
\section{Introduction}
\label{S:Intro}

The theory of domination in graphs and its variants has been an evergreen topic of research in graph theory over the past few decades.  Two vertices in a graph $G$ are \emph{neighbors} if they are adjacent in $G$. The \emph{open neighborhood} $N_G(v)$ of a vertex $v$ in $G$ is the set of neighbors of $v$, and its \emph{closed neighborhood} is $N_G[v] = N_G(v) \cup \{v\}$.

A \emph{dominating set} of a graph $G$ is a set $S$ of vertices of $G$ such that every vertex not in $S$ has a neighbor in $S$, where two vertices are neighbors if they are adjacent. The \emph{domination number} of $G$, denoted by $\gamma(G)$, is the minimum cardinality of a dominating set. The \emph{domatic number}, $d(G)$, of $G$ is the maximum number of disjoint dominating sets in $G$.

A \emph{total dominating set} of a graph $G$ with no isolated vertex is a set $S$ of vertices such that every vertex in $G$ has a neighbor in~$S$. The total domination number, $\gt(G)$, of $G$ is the minimum cardinality of a total dominating set of $G$. The \emph{total domatic number}, $d_t(G)$, of $G$ is the maximum number of disjoint total dominating sets~\cite{CDH} in $G$. This can also be considered as a coloring of the vertices such that every vertex has a neighbor of every color (and has been called the \emph{coupon coloring problem}~\cite{CKTV15coupon}).

A \emph{restrained dominating set} of $G$ is a dominating set $S$ of $G$ with the additional property that every vertex outside $S$ has a neighbor in~$S$. The restrained domination number, $\gr(G)$, of $G$ is the minimum cardinality of a restrained dominating set of $G$. The \emph{restrained domatic number}, $d_r(G)$, of $G$ is the maximum number of disjoint restrained dominating sets in~$G$.

The notion of domination and its variations in graphs has been studied a great deal; a rough estimate says that it occurs in more than 3000 papers to date. We refer the reader to the two so-called domination books by Haynes, Hedetniemi, and Slater~\cite{hhs1, hhs2}. In this paper, we study a variant of domination called $k$-tuple restrained domination in graphs.

Let $k \ge 1$ be an integer and let $G$ be a graph with $\delta(G) \ge k-1$. Harary and Haynes~\cite{HaHa00} defined a $k$-\emph{tuple dominating set}, abbreviated kD-set, in $G$ to be a set $S$ of vertices in $G$ such that every vertex of $S$ has at least $k-1$ vertices in $S$ and every vertex outside $S$ has at least $k$ neighbors in $S$; that is, $|N_G[v] \cap S| \ge k$ for every vertex $v$ of $G$. The $k$-\emph{tuple domination number} $\kDom(G)$ of $G$ is the minimum cardinality of a kD-set in $G$. As remarked in~\cite{HaHa00} the $1$-tuple domination number is the well-studied domination number. Thus, $\gamma(G) = \oneD(G)$. The \emph{$k$-tuple domatic number}, $\kdom(G)$, of $G$ is the maximum number of disjoint kD-sets in~$G$. A kD-set with cardinality $\kDom(G)$ is called a $\kDom$-set of $G$.

Let $k\ge 1$ be an integer and let $G$ be a graph with $\delta(G) \ge k$. A subset $S \subseteq V$ is a $k$-\emph{tuple total dominating set}, abbreviated kTD-set, in $G$ if every vertex in $G$ has at least $k$ neighbors in $S$; that is, $|N_G(v) \cap S| \ge k$ for every vertex $v$ of $G$. The minimum cardinality of a kTD-set in $G$ is the $k$-\emph{tuple total domination number} of $G$, denoted by $\kTDom(G)$. As remarked in~\cite{HeKa09}, the $1$-tuple total domination number is the well-studied total domination number. Thus, $\gamma_t(G) = \oneT(G)$. A kTD-set with cardinality $\kTDom(G)$ is called a $\kTDom$-set of $G$. The \emph{$k$-tuple total domatic number}, $\ktdom(G)$, of $G$ is the maximum number of disjoint kTD-sets in~$G$. The concept of $k$-tuple total domination in graphs was first studied by the authors in~\cite{HeKa09}.

Let $k \ge 1$ be an integer and let $G$ be a graph with $\delta(G) \ge k$. A \emph{$k$-tuple total restrained dominating set}, abbreviated kTRD-set, of $G$ is a kTD-set in $G$ with the additional property that every vertex outside $S$ has at least $k$ neighbors outside $S$. The minimum cardinality of a kTRD-set in $G$ is the $k$-\emph{tuple total restrained domination number} of $G$, denoted by $\kTRDom(G)$. The \emph{$k$-tuple total restrained domatic number}, $\krtdom(G)$, of $G$ is the maximum number of disjoint kTRD-sets in~$G$. A kTRD-set with cardinality $\kTRDom(G)$ is called a $\kTRDom$-set of $G$. The concept of $k$-tuple total restrained domination in graphs was first studied by Kazemi in~\cite{Kaz}.

In this paper, we introduce and study two new concepts, namely the $k$-tuple restrained domination number and the $k$-tuple restrained domatic number.

\begin{defn}
\emph{Let $k \ge 1$ be an integer and let $G$ be a graph with $\delta(G) \ge k-1$. A \emph{$k$-tuple restrained dominating set}, abbreviated kRD-set, of $G$ is a kD-set in $G$ with the additional property that every vertex outside $S$ has at least $k$ neighbors outside $S$; that is, $|N_G[v] \cap S| \ge k$ for every vertex $v \in V(G)$ and $|N_G(v) \cap (V(G)-S)| \ge k$ for every vertex $v \in V(G) \setminus S$. The minimum cardinality of a kRD-set in $G$ is the $k$-\emph{tuple restrained domination number} of $G$, denoted by $\kRDom(G)$. A kRD-set with cardinality $\kRDom(G)$ is called a $\kRDom$-set of $G$.}
\end{defn}

\begin{defn}
\emph{Let $k \ge 1$ be an integer and let $G$ be a graph with $\delta(G) \ge k-1$. A \emph{$k$-tuple restrained domatic partition}, abbreviated kRD-partition, of $G$ is a partition of $V(G)$ into kRD-sets. The \emph{$k$-tuple restrained domatic number}, $\krdom(G)$, of $G$ is the maximum number of disjoint kRD-sets in~$G$. }
\end{defn}

We remark that the $1$-tuple restrained domination number is the well-studied restrained domination number. Thus, $\gamma_r(G) = \oneR(G)$.

\subsection{Notation}

For notation and graph terminology, we will typically follow \cite{MHAYbookTD}. Throughout this paper, all graphs will be considered undirected, simple and finite. Specifically, let $G$ be a graph with vertex set $V(G)$ and edge set $E(G)$, and of order $n = |V(G|$ and size $m = |E(G)|$. If the graph $G$ is clear from the context, we simply write $V$ and $E$ rather than $V(G)$ and $E(G)$, and we write $G = (V,E)$.

For a set $S \subseteq V$, its \emph{open neighborhood} is the set $N_G(S) = \bigcup_{v \in S} N_G(v)$, and its \emph{closed neighborhood} is the set $N_G[S] = N_G(S) \cup S$.  For a set of vertices $S \subseteq V$, the subgraph of $G$ induced by $S$ is denoted by $G[S]$. The subgraph obtained from $G$ by deleting all vertices in $S$ and all edges incident with vertices in $S$ is denoted by $G - S$. If $S = \{v\}$, we simply write $G - v$ rather than $G - S$. If $X$ and $Y$ are vertex disjoint subsets of $G$, then $[X,Y]$ denotes the set of edges joining $X$ and $Y$ in $G$.

We denote the \emph{degree} of a vertex $v$ in $G$ by $d_G(v)$. Thus, $d_G(v) = |N_G(v)|$. The minimum and maximum degree among the vertices of $G$ is denoted by $\delta(G)$ and $\Delta(G)$, respectively. We denote the \emph{complement} of $G$ by $\barG$. Further, if $v$ is a vertex in $G$, then we denote the corresponding vertex in $\barG$ by $\barv$.
We denote the \emph{path}, \emph{cycle}, and \emph{complete graph} on $n$ vertices by $P_n$, $C_n$, and $K_n$, respectively, while $K_{n_1,n_2,\ldots,n_p}$ denotes a \emph{complete $p$-partite graph} with partite sets of size\textcolor{red}{s} $n_1, n_2, \ldots, n_p$. We use the standard notation $[k] = \{1,\ldots,k\}$.

\section{Fundamental Properties of $\kRDom(G)$}
\label{S:property1}

In this section, we present some fundamental properties of the $k$-tuple restrained domination number of a graph. The following observation relates the parameters defined in the introductory Section~\ref{S:Intro}.

\begin{obser}
\label{obs1}
If $G$ is a graph with $\delta(G) \ge k \ge 1$, then the following hold. \\[-29pt]
\begin{enumerate}
\item $d_{\times k,t}^{\, r}(G) \le d_{\times k}^{\, r}(G) \le d_{\times k}(G)$.
\item $d_{\times k,t}^{\, r}(G) \le d_{\times k,t}(G) \le d_{\times k}(G)$.
\item $\gamma_{\times k}(G) \le \gamma_{\times k}^{r}(G) \le \gamma_{\times k,t}^{r}(G)$.
\item $\gamma_{\times k}(G) \le \gamma_{\times k,t}(G) \le \gamma_{\times k,t}^{r}(G)$.
\end{enumerate}
\end{obser}

\begin{obser}
\label{obs2}
If $G$ is a graph of order~$n$ with $\delta(G) \ge k-1 \ge 1$, then the following hold. \\[-29pt]
\begin{enumerate}
\item $k \le \gamma_{\times k}^{r}(G) \le n$.
\item Every vertex of degree at most~$2k-1$ belongs to every kRD-set of $G$.
\item If $\delta(G) \le 2k - 1$, then $d_{\times k}^{\, r}(G) = 1$.
\item If $\Delta(G) \le 2k - 1$, then $\gamma_{\times k}^{r}(G) = n$.
\item If $\gamma_{\times k}^{r}(G) < n$, then $\kRDom(G) \le n - k - 1$ and $n \ge 2k - 1$.
\end{enumerate}
\end{obser}


The $k$-\emph{join} $F \circ_{k} H$  of a graph $F$ to a graph $H$ of order at least $k$ is defined in~\cite{HeKa09} to be the graph obtained from the disjoint union of $F$ and $H$ by joining each vertex of $F$ to at least~$k$ vertices of $H$. The order of the graph $H$ in the $k$-join $G \circ_{k} H$ we call the \emph{size} of the $k$-join. Further if $\delta(F) \ge k$ and $\delta(H) \ge k-1$, then we call the $k$-join $F \circ_{k} H$ a \emph{good} $k$-\emph{join}. We say that a graph $G$ has a $k$-join if $G \cong F \circ_{k} H$ for some $k$-join $F \circ_{k} H$.
If $F \circ_{k} H$ is a $k$-join and every vertex of $F$ is joined to exactly $k$ vertices of $H$, then we denote the $k$-join by $F \circ_{*k} H$.

\begin{thm}
\label{kDm}
Let $G$ be a graph with $\delta(G) \ge k-1 \ge 1$ and let $\ell$ be an integer such that~$k \le \ell < n$. Then, $\kRDom(G) = \ell$ if and only if $G$ has a good $k$-join and the minimum size of a good $k$-join in $G$ is~$\ell$.
\end{thm}
\proof Suppose that $\kRDom(G) = \ell$. Let $S$ be a $\kRDom$-set of $G$, and so $S$ is a kRD-set of $G$ and $|S| = \ell$. By supposition, $\ell < n$, and so $S$ is a proper subset of $V(G)$. Letting $H = G[S]$ and $F = G[V(G) \setminus S]$, we note that $F \circ_{k} H$ is a good $k$-join of size $n(H) = |S| = \ell$ and $G \cong F \circ_{k} H$. If $F' \circ_{k} H'$ is an arbitrary good $k$-join in $G$, then $V(H')$ is a $\kRDom$-set of $G$, implying that $\ell = \kRDom(G) \le n(H')$. Thus, the minimum size of a good $k$-join in $G$ is~$\ell$.

Conversely, suppose that $G$ has a good $k$-join and the minimum size of a good $k$-join in $G$ is~$\ell$. Thus, $G \cong F \circ_{k} H$ for some good $k$-join $F \circ_{k} H$ where $n(H) = \ell$. The set $V(H)$ is a $\kRDom$-set of $G$, implying that $\kRDom(G) \le \ell$.  If $\kRDom(G) < \ell$, then let $S'$ be a kRD-set of $G$, and so $|S| < \ell$. In this case, letting $H' = G[S]$ and $F' = G[V(G) \setminus S]$, we note that $F' \circ_{k} H'$ is a good $k$-join of size $n(H') = |S| < \ell$ and $G \cong F' \circ_{k} H'$, contradicting the fact that the minimum size of a good $k$-join in $G$ is~$\ell$. Hence, $\kRDom(G) = \ell$.~\qed

\medskip
As a consequence of Theorem~\ref{kDm}, we have the following result.

\begin{cor}
\label{c:kD}
Let $G$ be a graph with $\delta(G) \ge k-1 \ge 1$.  Then, $\kRDom(G) = k$ if and only if $G \cong K_{k}$ or $G \cong F \circ_{k} K_{k}$ for some graph $F$ with $\delta(F) \ge k$.
\end{cor}

Next we present a lower bound on the $k$-tuple restrained domination number of a graph in terms of its order and size.

\begin{thm}
\label{Lower1}
If $G$ is a graph of order $n$ and size $m$ with $\delta(G) \ge k-1$, then
\begin{equation}
\label{3.1}
\kRDom(G) \ge \frac{3kn-2m}{2k+1},
\end{equation}
with equality if and only $G \cong F \circ_{*k} H$ for some good $k$-join $F \circ_{*k} H$ where $H$ is a \emph{(}$k-1$\emph{)}-regular graph of order~$\kRDom(G)$ and $F$ is a $k$-regular graph of order $n - \kRDom(G)$.
\end{thm}
\proof
Let $S$ be a $\kRDom$-set of $G$, and so $S$ is a kRD-set of $G$ and $|S| = \kRDom(G)$. Let $\barS = V(G) \setminus S$, and so $|\barS| = n-\gamma_{\times k}^{r}(G)$. Let $H = G[S]$ and $F = G[\barS]$. Since $S$ is a kRD-set of $G$, we note that $\delta(H) \ge k-1$, $\delta(F) \ge k$ and each vertex in $\barS$ has at least~$k$ neighbors in $S$. Thus letting $m_1 = m(H)$, $m_2 = m(F)$ and $m_3 = |[S,\barS]|$, we have
\[
\begin{array}{lll}
m_1 & \ge & \frac{1}{2}(k-1)\kRDom(G), \1 \\
m_2 & \ge & \frac{1}{2} k(n-\kRDom(G)), \1 \\
m_3 & \ge & k(n-\kRDom(G)).
\end{array}
\]
Thus,
\[
m = m_1 + m_2 + m_3 \ge \frac{3}{2}kn - \frac{1}{2} (2k+1) \kRDom(G),
\]

\noindent
or, equivalently, $\kRDom(G) \ge (3kn-2m)/(2k+1)$. This establishes the desired lower bound. Suppose that we have equality in this lower bound. This implies that the above inequalities are all equalities; that is,
\[
\begin{array}{lll}
m_1 & = & \frac{1}{2}(k-1)\kRDom(G), \1 \\
m_2 & = & \frac{1}{2} k(n-\kRDom(G)), \1 \\
m_3 & = & k(n-\kRDom(G)).
\end{array}
\]

The first and second equalities implies that $H$ is a ($k-1$)-regular graph of order~$|S| = \kRDom(G)$ and $F$ is a $k$-regular graph of order $n - \kRDom(G)$, respectively, while the third equality implies that every vertex of $F$ is adjacent to exactly $k$ vertices of $H$, and so $G \cong F \circ_{*k} H$.~\qed

\medskip
Equality in the bound in Theorem~\ref{Lower1} is achieved, for example, by the complete graph $K_{2k+1}$ of order $2k+1$ which satisfies $\kRDom(K_{2k+1}) = k$.
Recall that $\gamma_r(G) = \oneR(G)$. Thus in the special case of Theorem~\ref{Lower1} when $k = 1$, we have the following lower bound on the restrained domination number of a graph.

\begin{cor}
\label{c:restDom}
If $G$ is a graph of order $n$ and size $m$, then $\gamma_r(G) \ge n- \frac{2}{3}m$.
\end{cor}

\section{Special Classes of Graphs}
\label{S:class}

In this section, we determine the $k$-tuple restrained domination number of special classes of graphs. We first consider a complete graph $K_n$ on $n$ vertices. By Observation~\ref{obs2}(d), $\kRDom(K_{n}) = n$ for $n \le 2k$. For $n \ge 2k+1$, every $k$-element subset of vertices in $K_n$ is a kRD-set of $K_n$, and so $\kRDom(K_{n}) \le k$. By Observation~\ref{obs2}(a), $\kRDom(K_{n}) \ge k$. Consequently, in this case $\kRDom(K_{n}) = k$. We state the result formally as follows.

\begin{obser}
\label{o:Kn}
For integers $n \ge k \ge 1$,
\[
\kRDom(K_n)= \left\{
\begin{array}{cc}
n & \mbox{if }n\le 2k, \\
k & \mbox{otherwise.}
\end{array}
\right.
\]
\end{obser}

As shown in~\cite{DHHLM}, the restrained domination number of a cycle $C_n$ for all $n \ge 3$ is given by $\oneR(C_n) = \gamma_r(C_n) = n - 2 \lfloor n/3 \rfloor$. By Observation~\ref{obs2}(d), $\kRDom(C_n) = n$ for all $n \ge 3$. We next determine the restrained domination number of the complement $\barC_n$ of a cycle $C_n$. Since $\barC_4 \cong 2K_2$, by Observation~\ref{o:Kn} we note that $\gamma_r(\barC_4) = 4$. Since $\barC_5 \cong C_5$, we note that $\gamma_r(\barC_5) = \gamma_r(C_5) = 3$. Since $n \ge 6$, if we let $V(C_n) = \{i \mid i \in [n] \}$ and $E(C_{n})=\{ij \mid i \in [n] \mbox{ and } j \equiv i+1 \, (\modo \, n) \}$, then $\{\overline{1},\overline {4}\}$ is a restrained dominating set of $\barC_n$, and so $\gamma_r(\barC_n) \le 2$. Noting that in this case $\gamma_r(\barC_n) \ge \gamma(\barC_n) = 2$, we deduce that $\gamma_r(\barC_n) = 2$. We state the result formally as follows.

\begin{obser}
\label{o:Cnr}
For $n \ge 4$ an integer,
\[
\oneR(\barC_n) = \gamma_r(\barC_n) = \left\{
\begin{array}{ll}
4 & \mbox{if }n=4, \\
3 & \mbox{if }n=5, \\
2 & \mbox{otherwise. }%
\end{array}%
\right.
\]
\end{obser}

We next determine the $k$-tuple restrained domination number of the complement $\barC_n$ of a cycle $C_n$ for $k \ge 2$ and $n \ge 5$.
\begin{prop}
\label{kRCn-}
For integers $n \ge k + 3 \ge 5$,
\[
\kRDom(\barC_{n})=\left\{
\begin{array}{cl}
n & \mbox{if } \, n\le 2k+2, \\
k+1 & \mbox{if } \, n \ge 2k+3.
\end{array}%
\right.
\]
\end{prop}
\proof
Since $\barC_n$ is $(n-3)$-regular, Observation~\ref{obs2}(d) implies that $\kRDom(\barC_n) = n$ for $n \le 2k+2$. Hence we may assume that $n \ge 2k+3$, for otherwise the desired result follows. By Observation~\ref{obs2}(a), $\kRDom(\barC_n) \ge k$. Let $S$ be a $\kRDom$-set of $\barC_n$, and so $S$ is a kRD-set of $\barC_n$ and $|S| = \kRDom(\barC_n)$. If $|S| = k$, then $G[S] \cong K_k$ and every vertex in $V(G) \setminus S$ is adjacent to every vertex in $S$, implying that $d_{\barC_n}(v) = n-1$ for all $v \in S$, contradicting the fact that $\barC_n$ is $(n-3)$-regular. Hence, $\kRDom(\barC_n) = |S| \ge k + 1$. Recall that  $n \ge 2k+3$. Let
\[
S^k_{\odd} = \bigcup_{i = 0}^{\lfloor \frac{k}{2} \rfloor} \{\, \overline{4i+1}, \overline{4i+2} \, \} \hspace*{0.5cm} \mbox{and} \hspace*{0.5cm} S^k_{\even} = \{\, \overline{1} \, \} \cup \left( \, \bigcup_{i = 1}^{\lfloor \frac{k}{2} \rfloor} \{\, \overline{4i}, \overline{4i+1} \, \} \right).
\]

For $k \ge 3$ odd, the set $S^k_{\odd}$ is a kRD-set of $\barC_n$, and so $\kRDom(\barC_n) \le |S^k_{\odd}| = k+1$. For $k \ge 2$ even, the set $S^k_{\even}$ is a kRD-set of $\barC_n$, and so $\kRDom(\barC_n) \le |S^k_{\even}| = k+1$. In both cases, $\kRDom(\barC_n) \le k+1$. Consequently, $\kRDom(\barC_n) = k + 1$.~\qed

\medskip
We next consider the $k$-tuple restrained domination number of a bipartite graph.

\begin{thm}
\label{p:bip1}
Let $G$ be a bipartite graph of order $n$ with $\delta(G)\ge k-1\ge 1$. Further, let $X$ and $Y$ be the partite sets of $G$, and let $\Delta_X$ and $\Delta_Y$ be the maximum degree among the vertices of $X$ and $Y$, respectively, in $G$. Then the following holds.
\\[-24pt]
\begin{enumerate}
\item If $\min \{ \Delta_X, \Delta_Y\} \le 2k-1$, then $\kRDom(G) = n$. \1
\item $\kRDom(G) \ge 2k-2$, with equality if and only if $G \cong K_{k-1,k-1}$. \1
\item If $|X| = k-1$, then $G \cong K_{k-1,|Y|}$. \1
\item $\kRDom(G) = 2k-1$ if and only if $G \cong K_{k-1,k}$.
\end{enumerate}
\end{thm}
\proof (a) Suppose that $\min \{ \Delta_X, \Delta_Y\} \le 2k-1$. Renaming the partite sets if necessary, we may assume that $\Delta_X = \min \{ \Delta_X, \Delta_Y\}$, and so $\Delta_X \le 2k-1$. Let $S$ be a $\kRDom$-set of $G$. If $X$ contains a vertex $v$ not is $S$, then $v$ would have at least $k$ neighbor in $Y$ that belong to the set $S$ and at least $k$ neighbors in $Y$ that do not belong to the set $S$, implying that $\Delta_X \ge d_G(v) \ge 2k$, a contradiction. Hence, $X \subseteq S$. This implies that every vertex in $Y$ has all its neighbors in the set $S$, and therefore every vertex of $Y$ belongs to the set $S$. Thus, $S = X \cup Y = V(G)$, and so $\kRDom(G) = |S| = n$.

(b) Let $S$ be a $\kRDom$-set of $G$, and let $x$ and $y$ be arbitrary vertices of $X$ and $Y$, respectively. Since $S$ is a kRD-set of $G$, we note that $|S \cap N(x)| \ge k-1$ and $|S \cap N(y)| \ge k-1$. Thus since $N(x) \cap N(y) = \emptyset$, we obtain $\kRDom(G) = |S| \ge |S \cap N(x)| + |S \cap N(y)| \ge 2k-2$. If $\kRDom(G) = 2k-2$, then we must have equality throughout this inequality chain, implying that $|S| = 2k-2$ and $|S \cap N(x)| = |S \cap N(y)| = k-1$. This in turn implies that $x \in S$ and $y \in S$. Since $x$ and $y$ are arbitrary vertices of $X$ and $Y$, respectively, we deduce that $X \subset S$ and $Y \subset S$ and therefore that $S = V(G) = X \cup Y$. Further, $|X| = |S \cap X| = k-1$ and $|Y| = |S \cap Y| = k-1$. Therefore, $G \cong K_{k-1,k-1}$. Conversely if $G \cong K_{k-1,k-1}$, then Observation~\ref{obs2}(d) implies that $\kRDom(G) = n = 2k-2$.

(c) Suppose that $|X| = k-1$. Let $S$ be a $\kRDom$-set of $G$, and let $y$ be arbitrary vertex of $Y$. Since $S$ is a kRD-set of $G$ and $S \cap N(y) \subseteq X$, we note that $k-1 = |X| \ge |S \cap N(y)| \ge k-1$, implying that $|X| = |S \cap N(y)| = k-1$ and $X \subseteq S$. This in turn implies that $y \in S$ and $y$ is adjacent to every vertex of $X$ in $G$, and so $G \cong K_{k-1,|Y|}$. This completes the proof of Part~(c). Part~(d) follows readily from Part~(c).~\qed

\medskip
As an immediate consequence of Theorem~\ref{p:bip1}(a), we have the following result.

\begin{cor}
\label{c:bip}
For integers $n \ge m \ge k-1 \ge 1$,
\[
\gamma_{\times k}^{r}(K_{n,m})=\left\{
\begin{array}{cl}
2k  & \mbox{if } m\ge 2k, \\
n+m & \mbox{if } m < 2k.
\end{array}%
\right.
\]
\end{cor}

Next we consider the $k$-tuple restrained domination number of a complete multipartite graph with at least three partite sets. For this purpose, we introduce the following notation. Let $G$ be a complete $p$-partite graph for some $p \ge 3$ and let $S$ be a $\kRDom$-set of $G$. We say that a partite set $X$ of $G$ is $S$-\emph{full} if every vertex in $X$ belongs to the set $S$; that is, if $X \subseteq S$. Let $f_{_S}(G)$ be the number of $S$-full partite sets in $G$. We note that if all $p$ partite sets of $G$ are $S$-full, then $f_{_S}(G) = p$. Let
\[
f(G) = \min \{p - f_{_S}(G) \mid S \mbox{ is a $\kRDom$-set of $G$ } \}.
\]

We remark that $f(G) = 0$ if and only if $\kRDom(G) = n(G)$. Moreover if a partite set $X$ of $G$ is not $S$-full for some $\kRDom$-set $S$ of $G$, then each vertex in $X \setminus S$ has at least~$k$ neighbors that not do not belong to $S$. Since these neighbors belong to partite sets different from $X$, this implies that at least one partite set of $G$ different from $X$ cannot be $S$-full. Thus if $f(G) > 0$, then $f(G) \ge 2$. We are now in a position to prove the following result.

\begin{thm}
\label{p:cpartite}
For $p \ge 3$, if $G$ is a complete $p$-partite graph of order $n$, then the following holds. \\[-20pt]
\begin{enumerate}
\item $\kRDom(G) \ge \left\lceil \frac{p(k-1)}{p-1} \right\rceil$.
\item If $\kRDom(G) < n$, then  $\kRDom(G) \le n - k - \left\lceil \frac{k}{f(G)-1} \right\rceil$.
\end{enumerate}
\end{thm}
\proof Let the complete $p$-partite graph $G$ have partite sets $X_1, X_2, \ldots, X_p$ where $|X_i| = n_i$ for $i \in [p]$. Thus, $G = K(n_{1},n_2, \cdots,n_{p})$. Let $S$ be a $\kRDom$-set of $G$, and let $\barS = V(G) \setminus S$. Further let $S_i = S \cap X_i$ and $\barS_i = X_i \setminus S$, and let $s_i = |S_i|$ for $i \in [p]$. Since every vertex in $X_i$ is adjacent to at least $k-1$ vertices in $S$ and these vertices all belong to $S \setminus X_i$, we note that
\[
|S| - s_i = \left( \sum_{j=1}^{p} s_{j} \right) - s_{i} \ge k-1
\]
for every $i \in [p]$. Thus,
\[
p|S| = \sum_{i=1}^{p} |S| \ge \sum_{i=1}^{p} (k + s_i - 1) = p(k-1) + |S|,
\]
or, equivalently, $|S| \ge p(k-1)/(p-1)$. This establishes Part~(a).

To prove Part~(b), suppose that $\kRDom(G) < n$. Thus, $f(G) > 0$, implying that $f(G) \ge 2$. Let $S$ be a $\kRDom$-set of $G$ such that $f(G) = p - f_{_S}(G)$. Renaming the partite sets $X_1, X_2, \ldots, X_p$, if necessary, we may assume that $s_i < n_i$ for $i \in [f(G)]$. Thus, each partite sets $X_i$ contains a vertex that belongs to $\barS$ for $i \in  [f(G)]$, while each partite sets $X_j$ is $S$-full for $j \in \{f(G) + 1, \ldots, p\}$. For each $j \in  [f(G)]$, let $x_j$ be a vertex of $X_j$ that belongs to $\barS_j$. Since $S$ is a kRD-set of $G$, the vertex $w_j$ has at least $k$ neighbors that belong to the set $\barS$. Since each such neighbor of $w_j$ belong to the set $\barS_i$ for some $i \in  [f(G)] \setminus \{j\}$, and since $N(w_{j}) \cap \barS_i = \barS_i$, we note that
\[
\begin{array}{lll}
k & \le & |N(w_{j}) \cap \barS| \2 \\
& = & \sum\limits_{i=1,i \ne j}^{f(G)} |N(w_{j}) \cap \barS_i| \2 \\
& = & \sum\limits_{i=1,i \ne j}^{f(G)} |\barS_i| \2 \\
& = & |\barS| - |\barS_j|.
\end{array}
\]
Thus,
\[
\begin{array}{lll}
k \cdot f(G) = \sum\limits_{j=1}^{f(G)} k  & \le & \sum\limits_{j=1}^{f(G)} (|\barS| - |\barS_j|) \2 \\
& = & f(G)\cdot |\barS| - |\barS| \1 \\
& = & (f(G) - 1) |\barS| \1 \\
& = & (f(G) - 1)(n - |S|),
\end{array}
\]
implying that
\[
\kRDom(G) = |S| \le n - k - \frac{k}{f(G) - 1},
\]
Since $\kRDom(G)$ is an integer, we therefore have that
\[
\kRDom(G) \le n - k - \left\lceil \frac{k}{f(G) - 1} \right\rceil,
\]
which completes the proof of Part~(b).~\qed

\section{Bounds on $\krdom(G)$}
\label{S:property2}

In this section, we present some fundamental properties of the $k$-tuple restrained domatic number of a graph. We first determine the $k$-tuple restrained domatic number of a complete graph.

\begin{obser}
\label{obsKn}
For $n \ge 2$ and $n \ge k \ge 1$, we have $\krdom(K_n) = \lfloor \frac{n}{k} \rfloor$.
\end{obser}

We establish next an upper bound on the product of the $k$-tuple restrained domination number and $k$-tuple restrained domatic number of a graph.

\begin{thm}
\label{t:product1}
If $G$ is a graph of order $n$ with $\delta(G) \ge k-1$, then
\[
\kRDom(G) \cdot \krdom(G) \le n.
\]
Moreover, if $\kRDom(G) \cdot \krdom(G) = n$, then there exists a partition of $V(G)$ into $\krdom(G)$ sets each of which is a $\kRDom$-set of $G$.
\end{thm}
\proof
Let $d = \kRDom(G)$, and so $d$ is the maximum number of disjoint kRD-sets in~$G$. Let $(V_1,\ldots, V_d)$ be a partition of $V(G)$ into kRD-sets in~$G$. Thus, each set $V_i$ is a kRD-set of~$G$ for $i \in [d]$, and so
\[
\krdom(G) \cdot \kRDom(G) = d \cdot \kRDom(G) = \sum_{i=1}^d \kRDom(G) \le \sum_{i=1}^d |V_i| = n.
\]

If $\kRDom(G) \cdot \krdom(G) = n$, then we must have equality throughout the above inequality chain, implying that $\kRDom(G) = |V_i|$ for all $i \in [d]$. Hence in this case, each set $V_i$ is a $\kRDom$-set of $G$.~\qed

\medskip
As a consequence of Corollary~\ref{c:kD} and Theorem~\ref{t:product1}, we have the following result.

\begin{cor}
\label{c:domr}
If $G$ is a graph of order $n$ with $\delta(G) \ge k-1 \ge 1$, then $\krdom \le \frac{n}{k}$ with equality if and only if $G \cong K_{k}$ or $G \cong F \circ_{k} K_{k}$ for some graph $F$ with $\delta(F) \ge k$.
\end{cor}

As a consequence of Theorem~\ref{p:bip1}, we have the following improvement on the upper bound of Corollary~\ref{c:domr} for the class of bipartite graphs.

\begin{cor}
\label{c:dombip}
If $G$ is a bipartite graph of order $n$ with $\delta(G) \ge k-1 \ge 1$, then $\krdom(G) \le \frac{n}{2k}$, unless $G \cong K_{k-1,k-1}$ or $G \cong K_{k-1,k}$, in which case $\krdom(G) = 1$.
\end{cor}

The following result establishes an upper bound on the $k$-tuple restrained domatic number of a graph in terms of its minimum degree.

\begin{thm}
\label{t:delta}
If $G$ is a graph with $\delta = \delta(G) \ge k-1 \geq 1$, then $\krdom(G) \le \frac{\delta + 1}{k}$.
\end{thm}
\proof Let $d = \kRDom(G)$ and let $(V_1,\ldots, V_d)$ be a partition of $V(G)$ into kRD-sets in~$G$. If $d = 1$, then the result is immediate since $d = 1 = \frac{k}{k} \le \frac{\delta + 1}{k}$. Hence we may assume that $d \ge 2$. let $v$ be a vertex of minimum degree in $G$, and so $d_G(v) = \delta$. Renaming the sets $V_1, V_2, \ldots, V_d$ if necessary, we may assume that $v \in V_k$. Thus, $|N_G(v) \cap V_k| \ge k-1$ and $|N_G(v) \cap V_i| \ge k$ for all $i \in [d-1]$. Thus,
\[
\begin{array}{lll}
\delta & = & |N_G(v)| \\
& = & \sum\limits_{i=1}^{d} |N_G(v) \cap V_i|  \1 \\
& \ge & (d-1)k + (k-1) \1 \\
& = & dk - 1,
\end{array}
\] 
and so $\krdom(G) = d \le (\delta + 1)/k$.~\qed

\medskip

We next obtain Nordhaus-Gaddum type results on the $k$-tuple restrained domatic number.

\begin{thm}
\label{t:NG}
If $G$ is a graph of order $n$ such that $\min \{\delta(G),\delta(\barG) \}\ge k-1 \geq 1$, then 
\[
\krdom(G) + \krdom(\barG) \le \frac{n+1}{k}. 
\]
Further if $\krdom(G) + \krdom(\barG) = \frac{n+1}{k}$, then the following holds. \\[-24pt]
\begin{enumerate}
\item  $\krdom(G) = \krdom(\barG) = \frac{n+1}{2k}$. \1
\item Both $G$ and $\barG$ are $(\frac{n-1}{2})$-regular graphs.
\end{enumerate}
\end{thm}
\proof Let $d = \krdom(G)$, $\delta = \delta(G)$ and $\Delta = \Delta(G)$, and let $\bard = \krdom(\barG)$ and $\bardelta = \delta(\barG)$. Applying Theorem~\ref{t:delta} to the graphs $G$ and $\barG$ we have $d \le (\delta + 1)/k$ and $\bard \le (\bardelta + 1)/k$. Thus since $\delta + \bardelta = n-1$, we note that $d + \bard \le (n+1)/k$. This establishes the desired upper bound. Suppose that $d + \bard = (n+1)/k$. This implies that $d = (\delta + 1)/k$ and $\bard = (\bardelta + 1)/k$. Without loss of generality, we may assume that $d \ge \bard$. Thus,
\begin{equation}
\label{Eq1}
d \cdot k - 1 \ge \bard \cdot k - 1 = \bardelta = \Delta \ge \delta = d \cdot k - 1.
\end{equation}

Hence we must have equality throughout inequality chain~(\ref{Eq1}), implying that $d = \bard$ and therefore that $d = (n+1)/(2k)$. Further, $\bardelta = \Delta = \delta$. This in turn implies that
\begin{equation}
\label{Eq2}
\bard \cdot k - 1 = d \cdot k - 1 = \delta = \barDelta \ge \bardelta = \bard \cdot k - 1.
\end{equation}
Hence we must have equality throughout inequality chain~(\ref{Eq2}), implying that $\Delta = \delta = \barDelta = \bardelta$. Thus since $2\delta = \delta + \bardelta = n-1$, both $G$ and $\barG$ are $(\frac{n-1}{2})$-regular graphs.~\qed

\begin{thm}
\label{t:NG2}
If $G$ is a graph of order $n$ with $\delta(G) \ge k-1 \ge 1$, then 
\[
\kRDom(G) + \krdom(G) \le n+1.
\]
\end{thm}
\proof Let $\gamma = \kRDom(G)$ and let $d = \krdom(G)$. By Theorem~\ref{t:product1}, 
$\gamma + d \le \frac{n}{d} + d$. By Corollary~\ref{c:domr}, we note that $1 \le d \le \frac{n}{k}$. Using the fact that the function $f(d) = \frac{n}{d} + d$ is decreasing
for $1\le d \le \sqrt n$ and increasing for $\sqrt n \le d \le n/2$, the function $f(d)$ attains its maximum value at one of its end points, namely at $d = 1$ or $d = n/k$, implying that $f(d) \le n+1$. Thus, $\gamma + d \le \frac{n}{d} + d \le n + 1$.~\qed

\medskip
We remark that the upper bound in Theorem~\ref{t:NG2} is achieved, for example, by all graphs $G$ satisfying $\kRDom(G) = n$.
The following result establishes a lower bound on the $k$-tuple restrained domatic number of a graph in terms of its order and minimum degree.

\begin{thm}
\label{t:dom3}
If $G$ is a graph of order $n$ with $\delta = \delta(G) \ge k-1 \ge 1$, then
\[
\krdom(G) \ge \left\lfloor \frac{n}{k(n - \delta)} \right\rfloor.
\]
\end{thm}
\proof If $n < k(n-\delta)$, then the result is immediate since in this case $\lfloor \frac{n}{k(n - \delta)} \rfloor \le 1 \le \krdom(G)$. Hence we may assume that $n \geq k(n-\delta)$. Thus, $n = pk(n-\delta)+r$ for some integers $p \ge 1$ and $0 \le r \le k(n-\delta)-1$. Let $S$ be an arbitrary subset of vertices of $G$ with $|S| \ge k(n-\delta)$ and let $\barS = V(G) \setminus S$. We note that $|S| \ge (n-\delta) + (k-1)(n-\delta) \ge n - \delta + k - 1$, implying that $|\barS| = n - |S| \le n - (n - \delta + k - 1) = \delta - k + 1$. If $v \in S$, then $|N(v) \cap S| \ge \delta - |\barS| \ge \delta - (\delta - k + 1) = k - 1$, while if $v \in \barS$, then $|N(v) \cap S| \ge \delta - (|\barS| - 1) \ge \delta - (\delta - k) = k$. Therefore, the set $S$ is a kRD-set of $G$. This is true for every set of subset of vertices of $G$ of size at least~$k(n-\delta)$. Thus letting $(S_1,S_2, \ldots, S_p)$ be a partition of $V(G)$ where $|S_i| = k(n-\delta)$ for $i \in [p-1]$ and letting $|S_p| = k(n-\delta)+r$, we produce a kRD-partition of $G$ into $p$ kRD-sets. Thus, $\krdom(G) \ge p = \lfloor \frac{n}{k(n - \delta)} \rfloor$.~\qed

\medskip
We remark that the lower bound of Theorem~\ref{t:dom3} is achieved, for example, when $G = K_n$ where $n \ge 2$ and $n \ge k$ as shown by Observation~\ref{obsKn}. We next present a sufficient condition for the $k$-tuple restrained domatic number of a graph be equal to its $k$-tuple domatic number. For this purpose, we first recall a result from~\cite{Kaz}.

\begin{prop}{\rm (\cite{Kaz})}
\label{TRD=TD} 
If $G$ is a graph with $\delta(G) \ge k \ge 1$, then $\krtdom(G) = \ktdom(G)$.
\end{prop}

\begin{thm}
\label{RD=D} 
If $G$ is a graph with $\delta(G) \ge k \ge 1$ and $(\kdom(G),\ktdom(G)) \ne (2,1)$, then $\krdom(G) = \kdom(G)$.
\end{thm}
\proof By Observation~\ref{obs1}(a), $\kdom(G) \ge \krdom(G) \ge \krtdom(G)$. If $\kdom(G) = 1$, then $\krdom(G) = 1$. If $\kdom(G) = 2$ and $\krdom(G) = 1$, then $\krtdom(G) = 1$, contradicting the assumption that $(\kdom(G),\ktdom(G)) \ne (2,1)$. Hence if $\kdom(G) = 2$, then $\krdom(G) = 2$. Therefore we may assume that $\kdom(G) \ge 3$, for otherwise $\krdom(G) = \kdom(G)$ as desired.

Let $d = \kdom(G)$ and let $(V_1,V_2,\ldots, V_d)$ be a partition of $V(G)$ into disjoint kD-sets in~$G$, where by assumption $d \ge 3$. Let $D$ be an arbitrary set in the partition $(V_1,V_2,\ldots, V_d)$, and so $D = V_i$ for some $i \in [d]$. We show that $D$ is a kRD-set of $G$. For notational simplicity, we may assume that $D = V_1$. Let $\barD = V(G) \setminus D$ and let $v$ be an arbitrary vertex in $\barD$. Renaming the sets $V_2, \ldots,V_d$ if necessary, we may assume that $v \in V_3$. Since $V_1$ and $V_2$ are kD-sets in $G$, there exist $k$-element subsets $D_v^1$ and $D_v^2$ such that $D_v^1 \subseteq N(v) \cap D_1$ and $D_v^2 \subseteq N(v) \cap D_2$. We note that $D_v^1 \subseteq D$ and $D_v^2 \subseteq \barD$. Hence, the vertex $v$ has at least $k$ neighbors in $D$ and at least $k$ neighbors in $\barD$. This is true for every vertex $v$ in $\barD$. Moreover since the set $D$ is a kD-set of $G$, every vertex in $D$ has at least $k-1$ neighbors that belong to $D$ in the graph $G$. Therefore, the set $D$ is a kRD-set of $G$. This is true for every set in the partition $(V_1,V_2,\ldots, V_d)$, yielding a partition of $V(G)$ into $d$ disjoint kRD-sets in~$G$. Hence, $\krdom(G) \ge d$. By Observation~\ref{obs1}(a), $d = \kdom(G) \ge \krdom(G) \ge d$, implying that $\krdom(G)= \kdom(G)$.~\qed

\medskip
Since $\kdom(K_{2k+1})=2$ and $\ktdom(K_{2k+1}) = \kdom(K_{2k+1})=1$, we remark that the condition $(\kdom(G),\ktdom(G)) \ne (2,1)$ in Theorem~\ref{RD=D} is necessary. In the special case when $k = 1$, Theorem~\ref{RD=D} yields the following result.

\begin{cor}
\label{Zel1} 
If $G$ is a graph with $\delta(G) \ge 1$ and $(d(G),d_t(G)) \ne (2,1)$, then $d_r(G) = d(G)$.
\end{cor}

We close with a sufficient condition for $\krdom(G) = \kdom(G)$. For a graph $G$ with $\delta(G) \ge k-1 \ge 1$, we denote by $\kdoms(G)$ the maximum number of disjoint kD-sets in $G$ such that at least one set in the partition is a $\kDom$-set of $G$.

\begin{thm}
\label{kTR=kT} 
If $G$ is a graph with $\delta(G) \ge k - 1 \ge 1$ and $\kdoms(G) \ge 3$, then $\kRDom(G) = \kDom(G)$.
\end{thm}
\proof Let $d = \kdoms(G) \ge 3$ and let $(V_1,V_2,\ldots, V_d)$ be a partition of $V(G)$ into disjoint kD-sets in~$G$, where the set $V_1$ is a $\kDom$-set of $G$. Let $D = V_1$ and let $\barD = V(G) \setminus D$. Proceeding exactly as in the proof of Theorem~\ref{RD=D}, the set $D$ is a kRD-set of $G$, implying that $\kRDom(G) \le |D| = \kDom(G)$. Conversely since every kRD-set of $G$ is also a kD-set of $G$, we have $\kRDom(G) \ge \kDom(G)$. Consequently, $\kRDom(G) = \kDom(G)$.~\qed

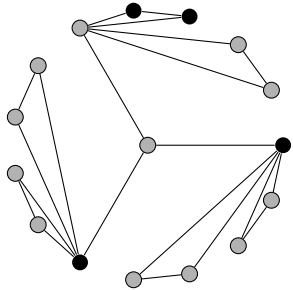
\begin{figure}[htp]
\centering
\begin{tikzpicture}[scale=1.8]
\tikzstyle{vertex}=[circle,fill=black,minimum size=6pt,inner sep=0pt]
\tikzstyle{vertex2}=[circle,draw,fill=black!30,minimum size=6pt,inner sep=0pt]

\node[vertex] (v1) at (0*360/5:1) {};
\node[vertex2] (v12a) at (0.33333*360/5:1) {};
\node[vertex2] (v12b) at (0.66666*360/5:1) {};
\node[vertex] (v2) at (1*360/5:1) {};
\node[vertex] (v23a) at (1.33333*360/5:1) {};
\node[vertex2] (v23b) at (1.66666*360/5:1) {};
\node[vertex2] (v3) at (2*360/5:1) {};
\node[vertex2] (v34a) at (2.33333*360/5:1) {};
\node[vertex2] (v34b) at (2.66666*360/5:1) {};
\node[vertex2] (v4) at (3*360/5:1) {};
\node[vertex] (v45a) at (3.33333*360/5:1) {};
\node[vertex2] (v45b) at (3.66666*360/5:1) {};
\node[vertex2] (v5) at (4*360/5:1) {};
\node[vertex2] (v15a) at (4.33333*360/5:1) {};
\node[vertex2] (v15b) at (4.66666*360/5:1) {};

\node[vertex2] (u1a) at (-0.0*0/5:0.0) {};

\draw  (v12b) -- (v23b) -- (v12a) -- (v12b) ;
\draw  (v23a) -- (v23b) -- (v2) -- (v23a) ;

\draw  (v4) -- (v45a) -- (v34b) -- (v4) ;
\draw  (v34a) -- (v45a) -- (v3) -- (v34a) ;

\draw  (v1) -- (v15b) -- (v15a) -- (v1) ;
\draw  (v1) -- (v5) -- (v45b) -- (v1) ;

\draw (v1) -- (u1a);
\draw (v45a) -- (u1a);
\draw (v23b) -- (u1a);

\end{tikzpicture}
\label{fi:TH10}
\caption{A graph $G$ with $\gamma(G)=3$, $\gamma_r(G)=4$ and $d(G) \ge 3$} 
\end{figure}

\medskip
For integers $k$ and $n$ where $5 \le 2k+1 \le n \le 3k-1$, we note that $\kRDom(K_n) = \kDom(K_n) = k$ and $\kdoms(K_n) = 2$, and therefore the converse of Theorem~\ref{kTR=kT} does not hold. We also remark that the condition $\kdoms(G) \ge 3$ in the statement of  Theorem~\ref{kTR=kT} cannot be replaced by the condition $\kdom(G) \ge 3$, as may be seen by considering the graph $G$ illustrated in Figure~\ref{fi:TH10}. In this example, $\gamma(G)=3$ and the set of three vertices of degree~$5$ is the unique $\gamma$-set of $G$. Further, $\gamma_r(G)=4$ and the set of the dark vertices shown in Figure~\ref{fi:TH10} is an example of a $\gamma_r$-set of $G$. Therefore, $\gamma_r(G) > \gamma(G)$. However, $d(G) \ge 3$.

\end{document}